\documentclass[preprint,12pt]{elsarticle}

\usepackage{amssymb}
\usepackage{amsmath}

\usepackage{graphicx}
\usepackage{caption}
\usepackage{subcaption}
\usepackage{float} 
\usepackage{xcolor,colortbl}

\journal{Computer Physics Communications}

\begin{document}

\begin{frontmatter}



\title{Projection-based stabilization for high-order incompressible flow solvers}


\author[label_author1]{Antonio Blanco-Casares} 
\author[label_author1]{Vishal Kumar}
\author[label_author1]{Daniel Mira}
\author[label_author1]{Oriol Lehmkuhl}

\affiliation[label_author1]{organization={Barcelona Supercomputing Center (BSC)},
            city={Barcelona},
            country={Spain}}

\begin{abstract}

This work presents a novel stabilization strategy for the Galerkin formulation of the incompressible Navier–Stokes equations, developed to achieve high accuracy while ensuring convergence and compatibility with high-order elements on unstructured meshes.
The numerical algorithm employs a fractional step method with carefully defined boundary conditions to obtain a consistent pressure field, enabling high-order temporal accuracy.
The proposed stabilization is seamlessly integrated into the algorithm and shares the same underlying principle as the natural stabilization inherent in the fractional step method, both rely on the difference between the gradient operator and its projection.
The numerical dissipation associated to the stabilization term is found to diminish with increasing polynomial order of the elements.
Numerical test cases confirm the effectiveness of the method, demonstrating convergence under mesh refinement and increasing polynomial order.
\end{abstract}


\begin{highlights}
\item Propose a low-dissipative and efficient stabilization method embedded in a standard fractional step method to solve incompressible flows. 
\item Demonstrate that the proposed method converges under mesh refinement and increasing element order, confirming its suitability for high-order elements, including spectral ones.
\end{highlights}


\end{frontmatter}


\section{Introduction}

Numerical simulations for incompressible flows are fundamental to advancing our understanding and predictive capabilities in a wide range of physical systems. The incompressible Navier–Stokes equations are inherently nonlinear and feature complex interactions between pressure and velocity fields, making them challenging to solve reliably. Inaccurate or unstable discretizations can lead to unphysical results or failure to capture critical flow features such as vortices and boundary layers. The development of high-fidelity numerical methods is essential not only for capturing the correct physics but also for enabling computationally efficient simulations that scale to complex geometries and high Reynolds number ($Re$) regimes.

The Continuous Galerkin approach, characterized by globally continuous basis functions, provides a natural and well-established framework for discretizing the incompressible Navier–Stokes equations \cite{huerta2003finite}.
The use of high-order elements has gained significant attention for improving the accuracy per degree of freedom (DoF). In particular, spectral elements combine the geometric flexibility of finite elements with the exponential convergence properties of spectral methods. Unlike low-order discretizations, which require dense meshes to achieve comparable accuracy, high-order elements can represent the solution with fewer DoFs and capturing fine-scale variations within each element. Although they require more numerical operations per DoF, their compatibility with modern computational architectures can offset these costs and, in many cases, lead to improved overall efficiency. These advantages makes them especially attractive for large-scale simulations, offering the potential for reduced computational cost while maintaining high fidelity, an advantage in scientific and engineering applications where both accuracy and performance are critical.

For uniform meshes, the Galerkin formulation can achieve optimal convergence, but for nonuniform meshes it might degrade. This degradation becomes especially pronounced in convection-dominated problems, where spurious oscillations arises. These oscillations come from the fact that pure Galerkin method leads to symmetric mass matrices, and for solving the non-symmetric operator of the convective term, it is known to be unstable.\cite{huerta2003finite}. On nonuniform meshes, irregular element alignment with the flow direction and local variations in mesh size can amplify these oscillations, escalating the loss of accuracy and a convergence order.

Among the stabilization techniques, the ones that preserve the order of the Pure Galerkin formulation are denoted as High-order stabilizations, among the variety of these methods, it stands out the Local projection stabilization (LPS) methods \cite{kuzmin2024weno2}, those address the shortcoming of more traditional stabilization residual-based techniques, like Streamline Upwind Petrov-Galerkin (SUPG) and Galerkin Least-squares (GLS) \cite{huerta2003finite}. This method uses the projection of gradients and do not require evaluating the whole residual of the equation due to its consistency property, which in presence of the diffusion operator it involves a second-order derivatives in the residual \cite{becker2001finite}.
The High-order stabilization removes the spurious oscillations, however in cases with presence of strong gradients, like supersonic or reacting flows, Gibbs oscillations will manifest at the discontinuities. Those cases require in addition the use of a shock capturing, which generally suppress the over-undershoots by adding an artificial viscosity at the required locations \cite{kuzmin2020entropy,kuzmin2023weno}.

For solving the incompressible Navier–Stokes equations, the coupling of velocity and pressure in the momentum equation in combination with the incompressibility constraint poses a major challenge \cite{fehn2017_stability}.
Fractional step methods offer an efficient solution by decoupling the velocity and pressure fields. They typically solve a convection–diffusion-type equation for the velocity and a Poisson equation for the pressure. On Dirichlet boundaries, these methods require Neumann boundary conditions for the pressure \cite{fehn2017_stability}. An overview of the various classes of these schemes can be found in Ref. \cite{guermond2006_overview} and \cite{spencer2005spectral}.

It was proved by Codina \cite{Codina2001pressure}, that the fractional step already introduces a natural stabilization on the pressure equation that depends on the time step size. The use of small time steps reduces the stability, this has been shown in many works such in Ref. \cite{fehn2017_stability}.
In advection dominated flows (high $Re$), the natural stabilization in the pressure is not enough to mitigate instabilities in the momentum equation.
The use of different forms of the convection term might help \cite{charnyi2017_emac}, but still in many other cases an additional stabilization is required.
This stabilization technique should preserve the order of the original Pure Galerkin formulation and introduce low diffusion to the system, in order to be energy-conserving.
This is a critical aspect for simulating turbulent ﬂows using Direct Numerical Simulation (DNS) or Large Eddy Simulation (LES). Where excessive numerical diffusion, in DNS, may obscure the effects of molecular diffusion, while in LES, it can mask the influence of the subgrid-scale model \cite{sanderse2013_energy}.
A low-diffusion, energy-conserving numerical method ensures that all dissipation present in the simulation arises from physically modeled processes, not artificial numerical artifacts. This is essential to accurately capture the energy spectrum across all relevant scales.

In recent years, new families of stabilization techniques have emerged, and its application on the incompressible Navier-Stokes equations presents a promising area of study.
The objective of this paper is to develop a stable formulation that is well-suited for high-order elements, such as spectral elements. The proposed method should exhibit convergence with both increasing element order and mesh refinement, while at the same being computationally efficient and scalable, making it a viable option for practical large-scale simulations.

The paper is organized as follows. In the next section, the stabilization technique is presented for a generic conservation law. Then, we show the fractional step technique is for solving the coupling between pressure and velocity, and the proposed stabilization term is fitted into to the momentum equation. The formulation is validated through several benchmark problems, an inviscid 2D shear layer, the Taylor-Green vortex, and flow over an airfoil profile using a nonuniform mesh. For the airfoil case, results are compared against reference data from the literature.

\section{Stabilized FEM}
As the Navier-Stokes equations are a specific case of conservation laws, we first consider the general form. Let $\phi(\mathbf{x},t)$ be the scalar conserved quantity that depends on the dimensional space location $\mathbf{x} \in \mathbb{R}^d$, $d \in \{1,2,3\}$, and the time instant $t \geq 0$. An initial-boundary value problem can be stated as,

\begin{equation}
\left\{
\begin{aligned}
    \quad \frac{\partial \phi}{\partial t} + \nabla \cdot \mathbf{f}(\phi) &= 0 \quad &&\mathrm{in} \, \Omega, \\
    \phi(\mathbf{x} , 0) &= \phi_0 \quad &&\mathrm{in} \, \Omega, \\
    \phi(\mathbf{x} , t) &= \phi_D \quad &&\mathrm{on} \, \Gamma_D,
\end{aligned}
\right.
\label{eq:conserv_law}
\end{equation}
where $\mathbf{f(\phi)}$ is the flux function, $\phi_0$ the initial data and $\phi_D$ specifies Dirichlet conditions applied on the boundary portion $\Gamma_D \subset \partial \Omega$.
To disctretize Eq. \ref{eq:conserv_law} in the finite element space $V_h$ of polynomials of order $\mathsf{p} \in \mathbb{N}$, a given bounded domain $\Omega$ is discretized into a mesh $\Omega_h$ composed by elements $K_e$, $1 \leq e \leq n_e$.

Any function $v_h \in V_h$ can be represented as,
    $v_h = \sum^{n_h}_{i=1} \varphi_i v_i$,
where $n_h$ is the set of all the points of the mesh, and $\varphi_1 , . . . , \varphi_{n_h}$ are globally continuous basis functions associated with the mesh points. The quantity $\phi_h$, which is a discrete approximation of $\phi$, satisfies the boundary condition on $\Gamma_D$ with the precision of the mesh size $h$.
To derive the variational formulation, we multiply the governing equation by a test function $w_h \in V_h$, which vanish on $\Gamma_D$, and integrate over the domain $\Omega$, leading to a set of linearly independent algebraic equations that define the discrete problem,

\begin{equation}
    \sum^{n_e}_{e=1} \int_{K_e} w_h \left( \frac{\partial \phi_h}{\partial t} + \nabla \cdot \mathbf{f}(\phi_h) \right) \text{d}\mathbf{x} = 0, \, \quad \forall \, w_h \in V_h.
\end{equation}
The weak form ensures that the modeled equation is a first-order differential equation, if necessary moving part of the differentiation of $\phi_h$ to the weight function $w_h$ by means of the Divergence theorem.
Substituting $w_h = \{\varphi_1,\,...\,\varphi_{n_h}\}$, leads to the system,
\begin{equation}
    m_{ij} \frac{\partial \phi_j}{\partial t} + \sum^{n_e}_{e=1} \int_{K_e} \varphi_i \nabla \cdot \mathbf{f}(\phi_h) \,\text{d}\mathbf{x} = 0,
    \label{eq:spatial_discretized_form}
\end{equation}
where the consistent mass matrix is defined as,
\begin{equation}
    m_{ij} = \sum^{n_e}_{e=1}  \int_{K^e} \varphi_i \varphi_j \,\text{d}\mathbf{x}.
\end{equation}

The main difficulty is that Eq. \ref{eq:spatial_discretized_form}, obtained with Galerkin Finite Element Method, has a central difference nature, and it is well known that this scheme is unconditionally unstable when solved with convection operators, which are non-symmetric, as proven by stability analysis \cite{toro2013riemann}. In practice, solutions of convection dominated transport problems by the Galerkin method are often corrupted by spurious node-to-node oscillations. To achieve a proper convergence, we focus on methods based on adding a stabilization term ($s_h^e$),

\begin{equation}
    \sum^{n_e}_{e=1} \left[ \int_{K_e} w_h \left( \frac{\partial \phi_h}{\partial t} + \nabla \cdot \mathbf{f}(\phi_h) \right) \text{d}\mathbf{x}  + s^e_h(w_h, \phi_h) \right] = 0.
\end{equation}

\subsection{Stabilization methods}

An easy way to make the system conditionally stable is to transform the central difference to an upwind difference, that is equivalent to add a large amount of artificial diffusion. With this, a diffusion term is added to our semi-discretized expression and leads to,

\begin{equation}
    s^{e,\text{LO}}_L = \nu^{e,\text{FU}} \int_{K_e} \nabla w_h \cdot \nabla \phi_h \,d\mathbf{x},
\end{equation}
omitting the natural Neumann boundary condition that arise when deriving the weak form. The first order upwind viscosity is given by \cite{kuzmin2020entropy},
\begin{equation}
    \nu^{e,\text{FU}} = \frac{(h^e/\mathsf{p}) \, ||\mathbf{f}'||_{L^\infty(K^e)}}{2},
    \label{eq:fu_visc}
\end{equation}
in which $h_e$ is the element size, and $\mathbf{f}'$ is the derivative of the flux obtained in the given element, which represents the propagation speed of the information.

The problem with this method is that is very diffusive, and the order of the original Pure Galerkin formulation is reduced, for that reason we denote it as a Low-order method, as in Ref. \cite{kuzmin2023weno}.

Alternatively, high-order methods like LPS can preserves better strong gradients. In LPS, the difference between the gradient and its projection captures oscillations caused by solving the non-symmetric operator term with a symmetric mass matrix.
The stabilization term reads \cite{kuzmin2020entropy,kuzmin2024property},
\begin{equation}
    s^{e,\text{LPS}}_L = c_s \, \nu^{e,\text{FU}} \int_{K_e} \nabla w_h \cdot \left( \mathbf{g}_h(\phi_h) - \nabla \phi_h \right) \,\text{d}\mathbf{x},
    \label{eq:lps}
\end{equation}
where $c_s$ is a parameter to be tuned.
The term $\mathbf{g}_h(\phi_h)$ represents the $L^2(\Omega)$ projection of $\nabla \phi_h$ onto $V_h$, which is a continuous approximation to the gradient,

\begin{equation}
  m_{ij} \, \mathbf{g}_h(\phi_h) = \sum^{n_e}_{e=1} \int_{K_e} w_h \, \left( \nabla \phi_h \right) \, \text{d}\mathbf{x}.
\end{equation}

As a computationally inexpensive alternative to obtain the projection with standard Lagrange elements, Kuzmin and Quezada de Luna \cite{kuzmin2020entropy} employ the lumped mass matrix.

This stabilization method uses the difference between the discrete gradient and its continuous counterpart projected onto a smoother space. That difference represents the unresolved, small, or oscillatory scales that can lead to numerical instabilities.
Ref. \cite{john2006_vms} show the similarity of this method with the Variational multiscale method (VMS), in which they model the fluctuations of the unresolved scales using the projection operator.
By penalizing these unresolved components, we effectively filter the parts of the solution that are not well-represented in $V_h$.

In contrast to SUPG, the stabilization term of the LPS scheme adds diffusion not only along the streamlines but also in all crosswind directions \cite{lohmann2017flux}.
According to Lohmann et al. \cite{lohmann2017flux}, in the 1D pure advection problem with constant velocity, this method is equivalent to the SUPG with projection-approximated time derivative.
The use of the projection operator makes the stabilization term independent of the time integration method.

\subsection{Spectral finite elements}
\label{sec:spectral}

Spectral elements address certain shortcomings of classical high-order FEM. Instead of equally spaced nodes, it uses Gauss-Lobatto-Legendre (GLL) spacing to avoid Runge's phenomenon, which consists on large oscillations near the endpoints of the interpolation interval when increasing the polynomial order \cite{spencer2005spectral}.

The GLL spacing is not only used for the nodes, but also for the quadrature points, so the quadrature points coincide with the cell nodes and therefore the mass matrix is always diagonal. Although it does not achieve exact integration as the traditional Gauss-Legendre (GL) quadrature, convergence is still attained.

\begin{figure}[ht]
\centering
    \includegraphics[width=0.2\textwidth]{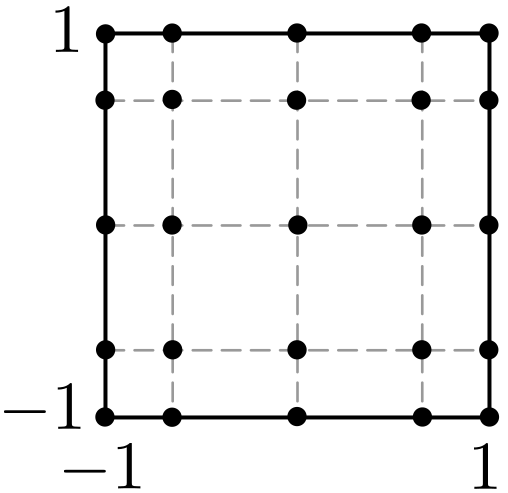}
    \caption{P4 quad element with GLL node spacing}
    \label{fig:p4_spectral}
\end{figure}

\section{Incompressible flows}
In incompressible flows, the continuity equation describes the divergence free property of the fluid ($\nabla \cdot \mathbf{u} = 0$) and momentum equation reads,
\begin{equation}
    \frac{\partial \mathbf{u}}{\partial t} + \mathcal{L}_N(\mathbf{u})  = \nabla \cdot \tau(\mathbf{u}) - \nabla p,
    \label{eq:mom}
\end{equation}
where we denote $p$ as the kinematic pressure \cite{huerta2003finite}. The viscous term is given by $\tau(\mathbf{u}) =  \nu \left( \nabla \mathbf{u} + (\nabla \mathbf{u})^\intercal \right)$, in which $\nu \geq 0$ is the kinematic viscosity.

Aliasing is a phenomenon that occurs when nonlinear terms in a PDE generate higher-frequency modes than what the numerical method can represent. This is the case of the nonlinear convection term ($\mathcal{L}_N(\mathbf{u})$), which requires degree $2\mathsf{p}-1$ but the numerical space can only represent up to degree $\mathsf{p}$. This error can lead to instabilities in under-resolved flows and is specially critical using high-order elements. It is well known that the conservative form is very sensitive to these errors, as an alternative the nonlinear convective term  can be expressed in different forms:
\begin{subequations}
\begin{align}
\text{Conservative:}\quad &\mathcal{L}_N = \nabla \cdot (\mathbf{u} \otimes \mathbf{u}) \\
\text{Non-conservative:}\quad &\mathcal{L}_N =  (\mathbf{u} \cdot \nabla) \mathbf{u} \\
\text{Skew-symmetric:}\quad &\mathcal{L}_N =
\frac{1}{2} \left( \nabla \cdot (\mathbf{u} \otimes \mathbf{u}) +(\mathbf{u} \cdot \nabla) \mathbf{u} \right) =
(\mathbf{u} \cdot \nabla) \mathbf{u} + \frac{1}{2} (\nabla \cdot \mathbf{u})\mathbf{u} \label{eq:skw}
\end{align}
\label{eq:convec}
\end{subequations}

The skew-symmetric form is defined as the combination of the conservative and non-conservative, and developing its expression it is equivalent to the non-conservative form plus an extra correction term that depends on the divergence of the velocity \cite{charnyi2017_emac}, which is theoretically zero by continuity but numerically might not.

\subsection{Time discretization of the Navier-Stokes equations}

The splitting between velocity and pressure from Eq. \ref{eq:mom} is performed by a non-incremental velocity-correction scheme \cite{karniadakis1991_splitting}, which according to Ref. \cite{guermond2006_overview} has similar convergence rates to its pressure-correction counterpart. It implies solving Equations \ref{eq:pred}-\ref{eq:diff} for a given time step size ($\delta t$),

\begin{align}
    \frac{\mathbf{u}^* - \mathbf{u}^n}{\delta t} &= - \mathcal{L}_N (\mathbf{u}^n) + \mathbf{s},
    \label{eq:pred} \\
    \nabla^2 p^{n+1} &= \frac{1}{\delta t} \nabla \cdot \mathbf{u}^*,
    \label{eq:poisson} \\
    \mathbf{u}^{**} &= \mathbf{u}^* - \delta t \, \nabla p^{n+1},
    \label{eq:corr} \\
    \mathbf{u}^{n+1} &= \mathbf{u}^{**} +  \nabla \cdot \tau(\mathbf{u}^{n+1}).
    \label{eq:diff}
\end{align}

This method, first predicts an intermediate velocity field $\mathbf{u}^*$ using the convective term. The pressure field is then obtained by solving a Poisson equation, which enforces incompressibility on the corrected velocity ($\nabla \cdot \mathbf{u}^{**} = 0$) and prescribes Dirichlet condition normal to the boundary \cite{karniadakis1991_splitting}, as will be described in the section below.
Unlike the the classical pressure-correction method, here once the incompressibility constrain has been applied, the last step (Eq. \ref{eq:diff}) solves the linear diffusion term with an implicit scheme. This mitigates potential inconsistencies near boundaries from the splitting, which might cause spurious pressure boundary layers on solid walls, and ensuring proper temporal accuracy.

The order of the temporal integration can be increased by means of multistage Runge–Kutta (RK) methods.
It is worth noting that we solve the Poisson equation for all the intermediate stages for the sake of consistency. A more efficient alternative would be the one proposed by Capuano et al. \cite{capuano2016approximate}, who developed an energy conserving RK-based approximated projection method in which the pressure correction is applied only at the RK final stage.

\subsection{Wall boundary condition}

The pressure in incompressible flows plays a major role, and imposing its proper boundary condition can leads to a high order in time scheme and minimize the effect of erroneous numerical pressure boundary layers induced by fractional step methods.
Karniadakis et al. \cite{karniadakis1991_splitting}, show an improved fractional step with the integration of the correct a Neumann BC, derived from evaluating the momentum equation in the normal of the domain boundary ($\mathbf{n}$), and they proved that it leads to a correct pressure solution.

The pressure condition at the wall can be derived from the momentum equation by applying the no-penetration constraint on the velocity,
\begin{equation}
    \left[-p \,\mathbf{I} + \nu \, (\nabla \mathbf{u} + (\nabla \mathbf{u})^\intercal) \right] \cdot \mathbf{n} = 0,
    \label{eq:zero_traction}
\end{equation}
noting that it only involves the viscous term. Ref. \cite{karniadakis1991_splitting} uses a special form of this condition,
\begin{equation}
    \nabla p^{n+1} \cdot \mathbf{n} = \frac{\partial p^{n+1}}{\partial n} = - \nu (\mathbf{n} \cdot \nabla \times (\nabla \times \mathbf{u}^{n+1})),
    \label{eq:p_bc}
\end{equation}
in which the operator $\nabla \times \nabla \times$ defines the so called rotational form of the boundary condition \cite{spencer2005spectral}. This form reinforces the incompressibility condition and yields unconditional stability \cite{guermond2006_overview}.

The Neumann boundary condition can be naturally incorporated into the Galerkin formulation. By integrating the second-order derivative in the Poisson equation (Eq. \ref{eq:poisson}) by parts, we obtain the weak form,

\begin{equation}
    \int_{\Omega_h} w_h \nabla^2 p_h^{n+1} \text{d}\mathbf{x} = - \int_{\Omega_h} \nabla w_h \nabla p_h^{n+1} \text{d}\mathbf{x} + \int_{\partial \Omega_h} w_h \left( \nabla p_h^{n+1} \cdot \mathbf{n} \right) \,\text{d}\mathbf{x},
\end{equation}

and its boundary term allows to directly impose the boundary condition specified in Eq. \ref{eq:p_bc}.
In alternative approaches, the pressure term is integrated by parts, canceling the Neumann boundary terms and naturally imposing the pressure boundary condition, which for this reason is sometimes denoted as ``do-nothing'' condition, like in Ref. \cite{braack2014_donothing}.

\subsection{Outflow boundary condition}

We assume  $\partial \Omega = \Gamma_D \cup \Gamma_0$, being $\Gamma_0$ section of the boundary in which neither the $\mathbf{u}$ or $p$, is known, denoted as the outflow boundary.
The zero traction condition, equivalent to Eq. \ref{eq:zero_traction}, applied on the outlet allows the fluid to flow freely out, as there is no stress-induced resistance.




The problem is that, in cases with strong backflow on $\Gamma_0$, sometimes caused by vortices at the outlet, the Neumann BC of based on the zero traction condition cannot provide a unique solution and it may become unstable or unphysical, because the velocity or pressure of the incoming fluid is not specified.

To address this problem, Brack and Mucha \cite{braack2014_donothing} and Dong et al. \cite{dong2014_outflow} propose similar approaches, both propose to introduce an additional term to Eq. \ref{eq:zero_traction} that compensates the possible inflows that may happen.
In this work, we adopt the method proposed by Dong et al., which includes a parameter to control the smoothing of the additional compensation term. The condition reads,
\begin{equation}
    -p \,\mathbf{n} + \nu \, (\nabla \mathbf{u} + (\nabla \mathbf{u})^\intercal)\cdot \mathbf{n} - \left(\frac{1}{2} |\mathbf{u}^{n+1}|^2 \,S_0(\mathbf{u})\right) \mathbf{n} = 0,
    \label{eq:dong}
\end{equation}

\begin{equation}
    S_0(\mathbf{u}) =  \frac{1}{2} \left(1- \mathrm{tanh}\left(\frac{\mathbf{n}\cdot \mathbf{u}}{U_0 \, \beta} \right) \right).
\end{equation}
The value $S_0$ is 1 if $\mathbf{n}\cdot \mathbf{u} <0$, and 0 otherwise. The hyperbolic tangent creates a smooth transition, which is controlled by the parameter $\beta$. $U_0$ is the characteristic velocity scale.


Dirichlet BC for the pressure in the outlet can be derived by multiplying Equation \ref{eq:dong} by the normal and considering $\mathbf{n} \cdot \mathbf{n} = 1$,

\begin{equation}
    p_h^{n+1} =  \nu \, \mathbf{n} \cdot \left( (\nabla \mathbf{u}^{n+1} + (\nabla \mathbf{u}^{n+1})^\intercal) \cdot \mathbf{n} \right) - \left(\frac{1}{2} |\mathbf{u}^{n+1}|^2 \, S_0(\mathbf{u}^{n+1})\right).
    \label{eq:pr_dirichlet}
\end{equation}
Since a Dirichlet condition is imposed on the pressure in $\Gamma_0$, the Neumann condition for the velocity is not applied on that part of the boundary.

\subsection{Stabilization of the fractional step method}

The work of Codina \cite{Codina2001pressure} proves that in the fractional step methods there is a natural pressure control which depends on the time step size.
In order to evidence it, we apply the incompressibility restrain (Eq. \ref{eq:corr}) in the FEM space,
\begin{equation}
    \mathbf{u}_h^{**} = \mathbf{u}_h^{*} - \delta t \, \mathbf{g}_h(p_h^{n+1}),
\end{equation}
where $\mathbf{g}_h$ is a continuous approximation of $\nabla p^{n+1}$. This expression can be used to substitute $\mathbf{u}_h^{*}$ in the Poisson equation (Eq. \ref{eq:poisson}), to obtain,
\begin{equation}
    \nabla \cdot \mathbf{u}_h^{**} = \delta t \, (\nabla^2 p_h^{n+1} - \nabla \cdot \mathbf{g}_h(p_h^{n+1})) = \delta t \,\nabla \cdot (\nabla p_h^{n+1} - \mathbf{g}_h(p_h^{n+1})).
    \label{eq:natural_stab}
\end{equation}

It is clear that divergence-free equation has an extra term added to the residual, which is the rhs of Eq. \ref{eq:natural_stab}.
The difference between the gradient and its projection represent the fluctuations \cite{john2016finite}, and the bigger they are, the more stabilizing effect, which will act as source term to the pressure.
Since the natural stabilization depends on $\delta t$, there is a lower bound for the time step for stability reasons.
The natural stabilization has the same order as the as the pressure interpolation error, which is $\mathcal{O}(\delta t \, h^{\mathsf{p}+1})$. In Ref. \cite{lehmkuhl2019low}, Lehmkuhl et al. prove that the error kinetic energy has the same error.

In many convection-dominated problems, the natural stabilization is not enough, and a proper stabilization method to the momentum equation is still required.
In this work we propose incorporating the high-order stabilization method to the prediction step, represented as $\mathbf{s}$ in Eq. \ref{eq:pred}.
Recalling the generic conservation law from Eq. \ref{eq:conserv_law}, in the momentum equation the conserved variable ($\phi$) is each component of the velocity vector. Therefore, the LPS stabilization term (Eq. \ref{eq:lps}) can be adjusted for each velocity component ($u_k \in \mathbf{u}$), becoming,
\begin{equation}
    s^{e,\text{LPS}}_{h,k} = c_s \nu^{e,\text{FU}} \int_{K_e} \nabla w_h \cdot ( \mathbf{g}_h(u_{h,k}) - \nabla u_{h,k}) \,\text{d}\mathbf{x}, \quad \nu^{e,\text{FU}} = \frac{(h^e/\mathsf{p}) ||\mathbf{u}_h||_{L^\infty(K^e)}}{2},
    \label{eq:inco_lps}
\end{equation}

It is interesting to observe the similarity between the natural stabilization from Eq. \ref{eq:natural_stab} and the LPS stabilization. Both are based on the residual of the gradient projection. If the natural stabilization based on the pressure gradient fails to deliver a smooth solution, the LPS will target the fluctuations of the velocity gradient.

Note that for incompressible flows, the derivative of the flux ($\mathbf{f'}$) in Eq. \ref{eq:fu_visc}, which is the propagation speed, is just the velocity of the fluid.
The additional parameter $c_s \in [0,1]$ can be used to adjust the level of high-order stabilization \cite{kuzmin2020entropy, kuzmin2023weno}.
In problems that employ shock capturing, this tool already introduces stabilizing effects throughout the domain, allowing a reduction in high-order dissipation.
This is not the case in incompressible flows and we set $c_s=1$ by default.

Using an additional projection-based stabilization method in incompressible flows is not something new, Codina and Blasco \cite{CodinaBlasco1997,codina2001implementation} and Becker and Braack \cite{becker2001finite} used the gradient of the discrete pressure. The novelty of this work relies on the adaptation of the high order stabilization described by Lohmann et al. \cite{lohmann2017flux} and Kuzmin and Quezada de Luna \cite{kuzmin2020entropy} for generic conservation laws and adapt it to incompressible flows using the velocity field.

\section{Test cases}

\subsection{Shear-layer problem}

To analyze kinetic energy conservation, we begin with an inviscid two-dimensional problem with periodic boundary conditions. The initial condition consists of a shear velocity profile centered in the domain and a sinusoidal velocity in the perpendicular direction, as described in Ref. \cite{sanderse2013_energy}.
From the start of the simulation the flow develops stables eddies, which end up resembling the 2D Taylor-Green vortex problem, and unlike the 3D version that will be shown below, the eddies do not decompose into smaller scales.

The computational domain spans $[0,2\pi]\times[0,2\pi]$,  with a uniform mesh of linear elements. The initial condition for $\mathbf{u}=(u,v)$ and $\mathbf{x}=(x,y)$ is defined as follows:
\begin{equation}
u = 
\begin{cases}
\tanh\left(\dfrac{y - \pi/2}{\delta}\right), & y \leq \pi, \\
\tanh\left(\dfrac{3\pi/2 - y}{\delta}\right), & y > \pi,
\end{cases}
\quad
v = \varepsilon \sin(x),
\end{equation}
in which $\delta =\pi/15$ and $\varepsilon =0.05$ \cite{sanderse2013_energy}.

\begin{figure}[ht]
    \captionsetup[subfigure]{justification=centering}
    \centering
    \begin{subfigure}[b]{0.3\textwidth}
        \centering
        \includegraphics[width=\textwidth]{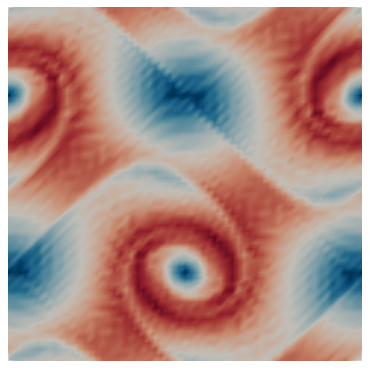}
        \caption{Pure Galerkin}
    \end{subfigure}
    \begin{subfigure}[b]{0.3\textwidth}
        \centering
        \includegraphics[width=\textwidth]{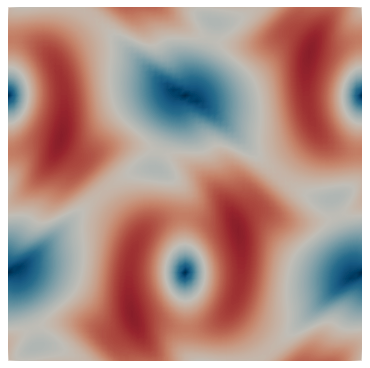}
        \caption{LPS}
    \end{subfigure}
    \begin{subfigure}[b]{0.06\textwidth}
        \centering
        \includegraphics[width=\textwidth]{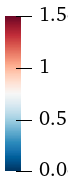}
    \end{subfigure}
    \caption{Shear layer problem at time $t=8$ solved with skew-symmetric formulation, P1 elements and a mesh $60 \times 60$ DoFs.}
    \label{fig:shear_layer}
\end{figure}

This case is solved using the different formulations of the convective term presented in Eq. \ref{eq:convec}, both with and without stabilization.
The velocity contours at time $t=8$ using skew-symmetric formulation are analyze in Figure \ref{fig:shear_layer}, capturing the moment right before stable eddies are fully developed. The time step is fixed to $\delta t = 5 \times 10^{-3}$.
The strong gradients observed, which might resemble discontinuities, are not aligned with the streamwise direction, instead, they arise from the imposed shear velocity. Therefore, the velocity gradient remains smooth along the streamlines,  which is consistent with the behavior expected in incompressible flows.
When applying LPS stabilization, the sharp gradients are noticeably smoothed, but at the same time, the small instabilities observed in the Pure Galerkin case are suppressed.

The average kinetic energy is computed as,
\begin{equation}
    E_k = \frac{1}{\Omega} \int_\Omega \frac{\mathbf{u}\cdot \mathbf{u}}{2}  \,\text{d}\mathbf{x}
    \label{eq:kinetic_energy}
\end{equation}

\begin{figure}[ht]
    \captionsetup[subfigure]{justification=centering}
    \centering
    \begin{subfigure}[b]{0.38\textwidth}
        \centering
        \includegraphics[width=\textwidth]{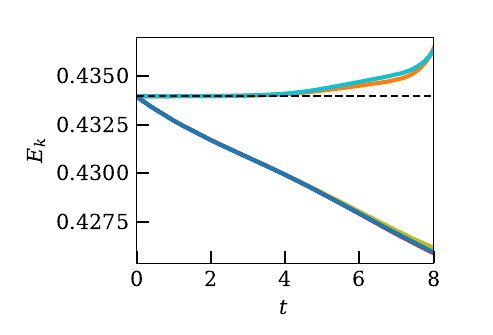}
        \caption{}
    \end{subfigure}
    \begin{subfigure}[b]{0.38\textwidth}
        \centering
        \includegraphics[width=\textwidth]{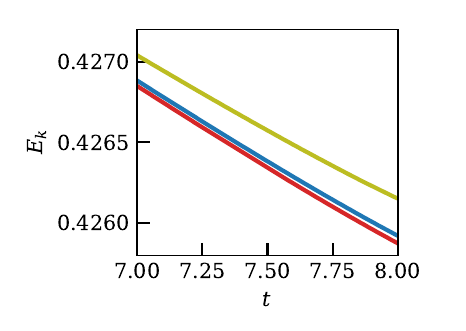}
        \caption{Zoom}
    \end{subfigure}
    \begin{subfigure}[b]{0.20\textwidth}
        \centering
        \includegraphics[width=\textwidth]{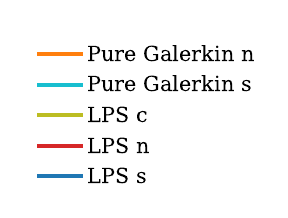}
    \end{subfigure}
    \caption{Kinetic energy evolution for the shear-layer problem using different convection formulations: conservative (c), non-conservative (n), and skew-symmetric (s).}
    \label{fig:shear-layer_Ek}
\end{figure}

Figure \ref{fig:shear-layer_Ek} shows the temporal evolution of kinetic energy for the shear-layer case.
The conservative formulation without stabilization is not included in the plot, as it leads to a rapid outbreak of instabilities that grow stronger than in the non-conservative or skew-symmetric cases, ultimately causing the simulation to diverge. However, the conservative formulation becomes stable when LPS stabilization is applied.

For the unstabilized cases, the kinetic energy remains approximately constant over time, but later increases as instabilities introduce numerical errors. In contrast, when stabilization is applied, the kinetic energy decreases due to the added numerical dissipation across all three convective formulations. As expected, the conservative formulation with stabilization preserves more energy than the non-conservative and skew-symmetric forms.

\subsection{3D Taylor-Green vortex}

\begin{figure}[ht]
    \captionsetup[subfigure]{justification=centering}
    \centering
    \begin{subfigure}[b]{0.3\textwidth}
        \centering
        \includegraphics[width=\textwidth]{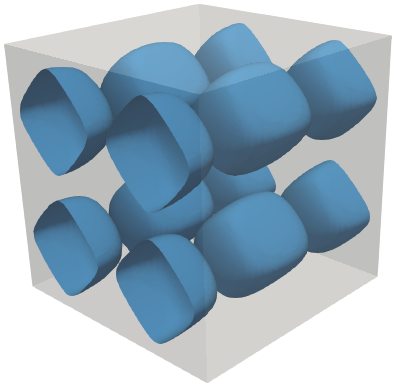}
        \caption{$t=0$}
    \end{subfigure}
    \begin{subfigure}[b]{0.3\textwidth}
        \centering
        \includegraphics[width=\textwidth]{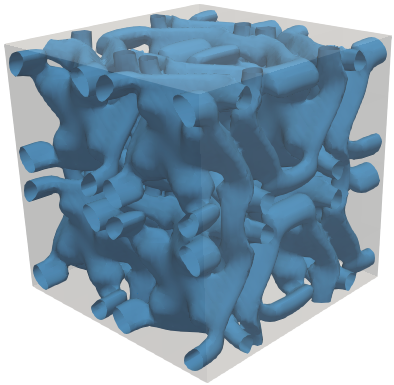}
        \caption{$t=5$}
    \end{subfigure}
    \begin{subfigure}[b]{0.3\textwidth}
        \centering
        \includegraphics[width=\textwidth]{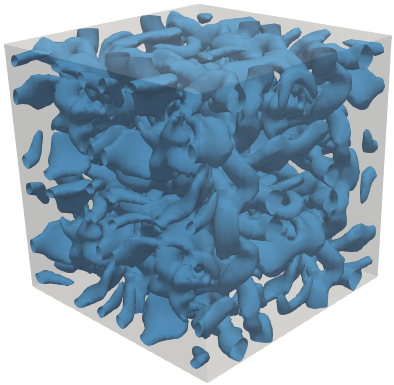}
        \caption{$t=10$}
    \end{subfigure}
    \caption{Q-criterion iso-surfaces of the Taylor-Green vortex problem using $64^3$ DoFs.}
    \label{fig:tgv_qcrit}
\end{figure}

The 3D Taylor-Green Vortex is a popular test case used to study turbulence and vortex dynamics, it serves as a benchmark problem for evaluating numerical methods, turbulence models, and the dissipation properties of high-fidelity simulations.
The domain consists on a cube with periodic boundary conditions. The Taylor-Green vortex is initialized with large-scale symmetric vortex structures that undergoes nonlinear interactions, leading to the breakdown of the vortices and the formation of turbulence.
Vortex stretching is a mechanism unique to three-dimensional flows, and is crucial in driving the transition from coherent vortices to fully developed turbulence \cite{lehmkuhl2019low}, as seen in Fig. \ref{fig:tgv_qcrit}.

The dimension of the cubic domain is $(2 \pi)^3$, and we set $Re = \frac{V_0 L_0}{\nu}= 1600$.
The flow is simulated for $t = 20\, t_c$, in which $t_c$ is the characteristic convective time $t_c = V_0/L_0$, in this case $L_0=1$ and $V_0=1$.
The initial condition is given by ($\mathbf{u}=(u,v,w)$ and $\mathbf{x}=(x,y,z)$) \cite{lehmkuhl2019low},
\begin{equation}
\begin{aligned}
u &= V_0 \sin(x) \cos(y) \cos(z), \\
v &= -V_0 \cos(x) \sin(y) \cos(z), \\
w &= 0.
\end{aligned}
\end{equation}

The problem is solved using three different uniform grids of $64^3$, $128^3$ and $256^3$ DoFs, using a set of standard Lagrange and spectral elements.

In addition to the integrated kinetic energy (Eq. \ref{eq:kinetic_energy}), we measure its dissipation rate ($\varepsilon$), which for incompressible flows can be computed from the integrated enstrophy ($\zeta$) \cite{debonis2013tgv},
\begin{equation}
    \varepsilon(\zeta) = 2 \nu \zeta, \quad \zeta = \frac{1}{\Omega} \int_\Omega \frac{\omega \cdot \omega}{2}  \,\text{d}\mathbf{x},
\end{equation}
and the vorticity is calculated from the velocity field, $\omega = \nabla \times \mathbf{u}$.

\begin{figure}[ht]
    \captionsetup[subfigure]{justification=centering}
    \centering
    \begin{subfigure}[b]{0.6\textwidth}
        \centering
        \includegraphics[width=\textwidth]{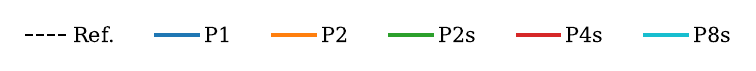}
    \end{subfigure}
    \begin{subfigure}[b]{0.62\textwidth}
        \centering
        \includegraphics[width=\textwidth]{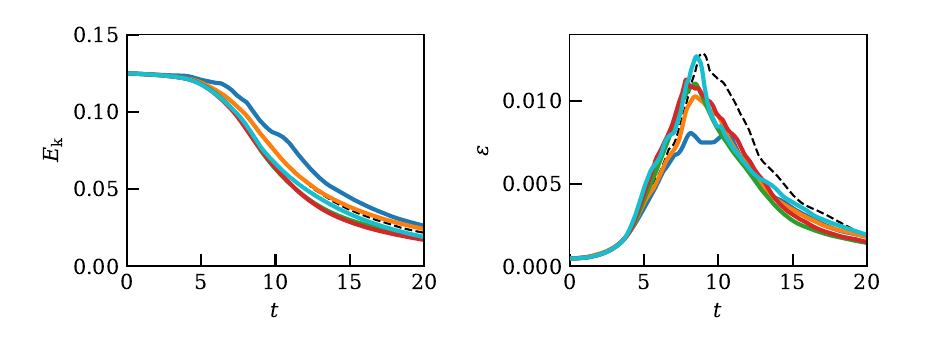}
        \caption{$64^3$ DoFs}
    \end{subfigure}
    \begin{subfigure}[b]{0.62\textwidth}
        \centering
        \includegraphics[width=\textwidth]{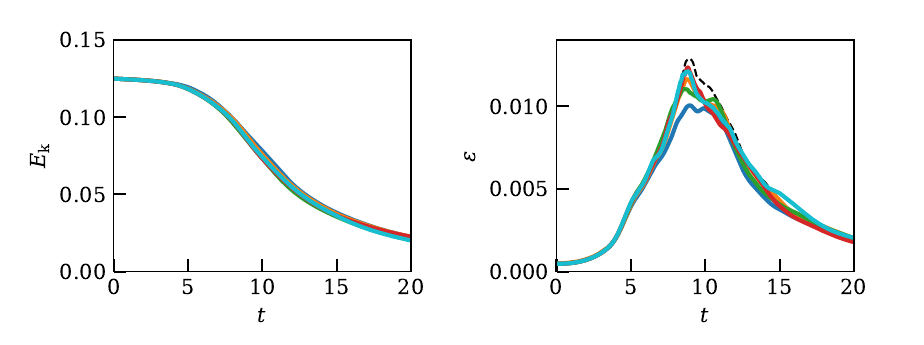}
        \caption{$128^3$ DoFs}
    \end{subfigure}
    \begin{subfigure}[b]{0.62\textwidth}
        \centering
        \includegraphics[width=\textwidth]{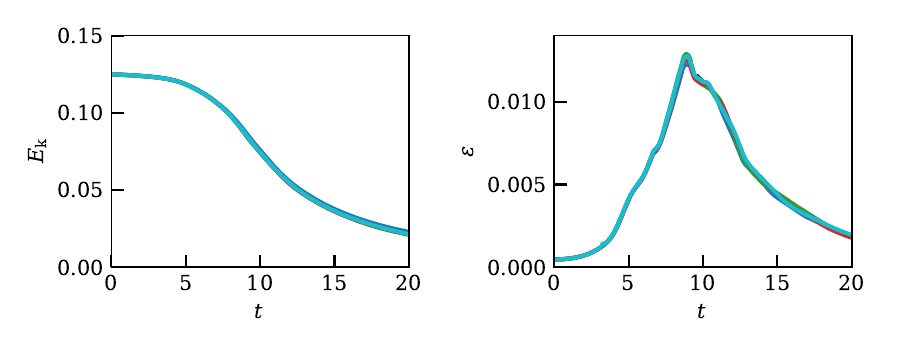}
        \caption{$256^3$ DoFs}
    \end{subfigure}

    \caption{Kinetic energy ($E_k$) and its dissipation rate ($\varepsilon$) in the Taylor-Green vortex problem using Pure Galerkin.}
    \label{fig:tgv_nolps}
\end{figure}

We first evaluate this test case with the Pure Galerkin method (Fig. \ref{fig:tgv_nolps}), noting that the physical viscosity already introduces a certain degree of stabilization.
It is clearly observed that in the $64^3$ DoFs mesh, the kinetic energy of some element types lies above the reference curve, provided in Ref. \cite{debonis2013tgv}. This behavior is attributed to the use of skew-symmetric formulation, if the conservative form of the convection term were used instead, the kinetic energy would consistently fall below the reference.
In the medium size mesh, the kinetic energy results appear to have converged; however, the dissipation rate, which is more sensitive, still shows noticeable deviations from the reference. In the finest mesh, the solutions for all element types have converged for both kinetic energy and dissipation rate.

\begin{figure}[ht]
    \captionsetup[subfigure]{justification=centering}
    \centering
    \begin{subfigure}[b]{0.6\textwidth}
        \centering
        \includegraphics[width=\textwidth]{figures/tgv/tgv_legend.pdf}
    \end{subfigure}
    \begin{subfigure}[b]{0.62\textwidth}
        \centering
        \includegraphics[width=\textwidth]{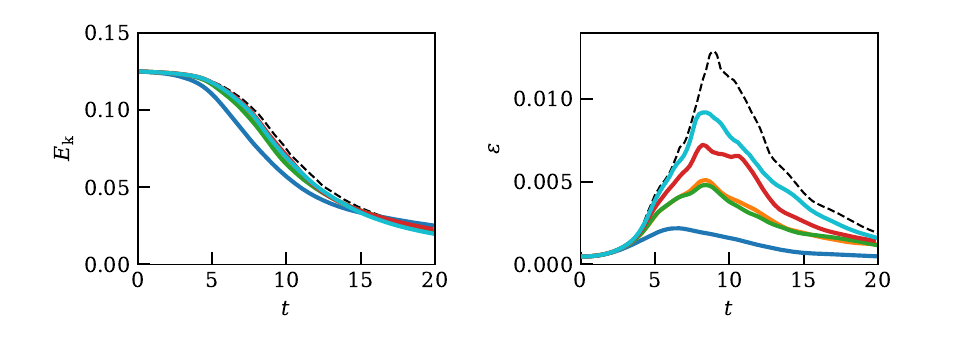}
        \caption{$64^3$ DoFs}
    \end{subfigure}
    \begin{subfigure}[b]{0.62\textwidth}
        \centering
        \includegraphics[width=\textwidth]{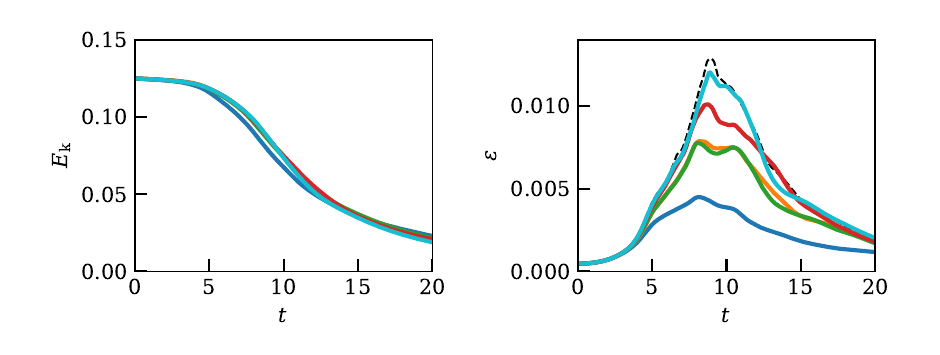}
        \caption{$128^3$ DoFs}
    \end{subfigure}
    \begin{subfigure}[b]{0.62\textwidth}
        \centering
        \includegraphics[width=\textwidth]{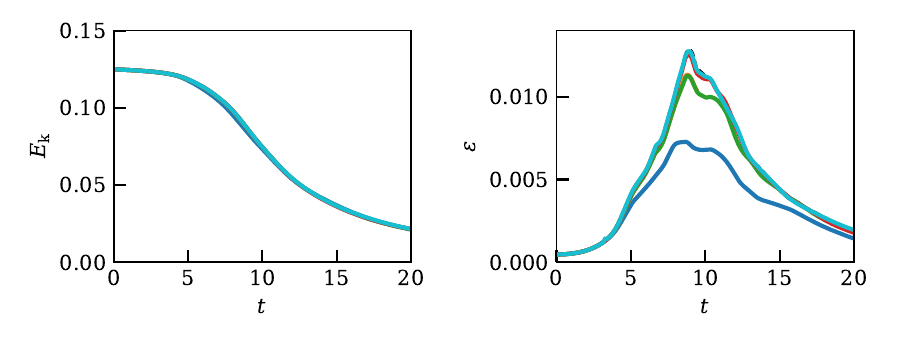}
        \caption{$256^3$ DoFs}
    \end{subfigure}
    \caption{Kinetic energy ($E_k$) and its dissipation rate ($\varepsilon$) in the Taylor-Green vortex problem using LPS stabilization.}
    \label{fig:tgv_lps}
\end{figure}

Figure \ref{fig:tgv_lps} presents the results applying the LPS stabilization method.
It is worth noting that the standard P2 elements match the spectral version P2s, this is a good indicator that both numerical approaches behave similarly.
Compared to the Pure Galerkin cases, a clear effect of the stabilization is the reduction of kinetic energy, and more significantly, the dissipation rate.

Interestingly, increasing the polynomial order reduces the dissipation introduced by the stabilization.
These results clearly demonstrate that the stabilized method converges to the reference solution with both increasing polynomial order and mesh refinement.
In fact, in the same way as the natural pressure stabilization, based on Codina and Blasco \cite{CodinaBlasco1997Afinite}, the error of the stabilization is $||\mathbf{g}_h(\phi_h) - \nabla \phi_h||_{L^2(\Omega)} = \mathcal{O}(h^{\mathsf{p}+1})$.
In P1, the gradient of the discrete approximation is a constant on each element, the difference term captures the part of the gradient that cannot be resolved by the continuous space. Since $\nabla u_h$ is more accurate when increasing the polynomial order, it introduces less distortion and the gradient difference will be lower, reducing the stabilization term and associated numerical diffusion, as proved by Fig. \ref{fig:tgv_lps}.

An alternative explanation is that when the quadrature points coincide with the nodes, the discrepancy between the discrete and continuous gradients vanishes at the interior nodes of the element. As the polynomial order increases, more points lie within the element, leading to a reduction of stabilization need for a fixed number of DoFs.



\subsection{Airfoil case}

For the final test case, we investigate the flow over NACA0012 at an angle of attack of $5.0^{\circ}$ and a Reynolds number of $Re = U_0 \,c / \nu= 5 \times 10^4$, being $U_0$ the free-stream velocity and $c$ the chord length.
The results are compared with the numerical solutions reported by Lehmkuhl et al. \cite{lehmkuhl2019low} and Zhang and Samtaney \cite{Zhang2016_airfoil} .

\begin{figure}[ht]
\centering
    \includegraphics[width=0.45\textwidth]{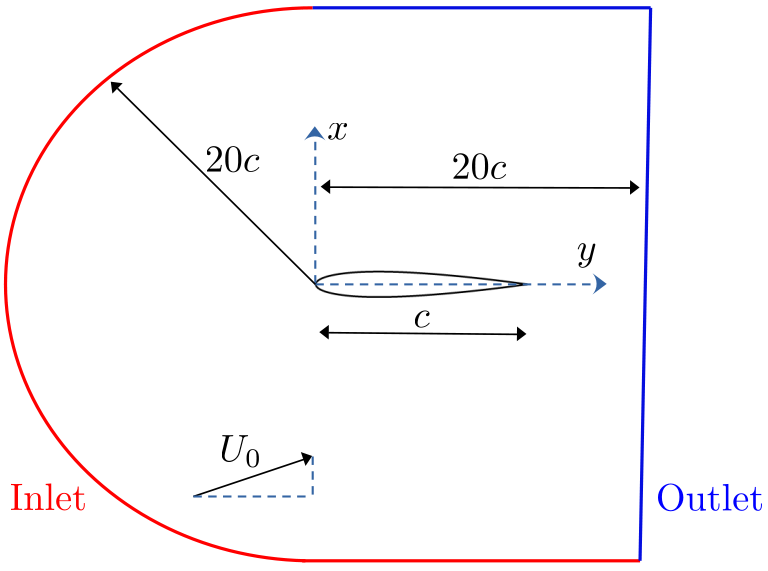}
    \caption{Computational domain of the Airfoil case}
    \label{fig:naca_domain}
\end{figure}

Following the setup in Ref. \cite{lehmkuhl2019low}, the computational domain consists on a three-dimensional volume with the profile illustrated in Fig. \ref{fig:naca_domain}. The domain has a spanwise depth of $0.2\,c$, with periodic boundary conditions in that direction. The given nonuniform mesh is composed of P4 spectral elements and has approximately 30 million DoFs.
The simulation employs a large-eddy simulation (LES) model to define the subgrid-scale viscosity. The solver uses skew-symmetric formulation of the convective term.

\begin{figure}[ht]
    \captionsetup[subfigure]{justification=centering}
    \centering
    \begin{subfigure}[b]{0.4\textwidth}
        \centering
        \includegraphics[width=\textwidth]{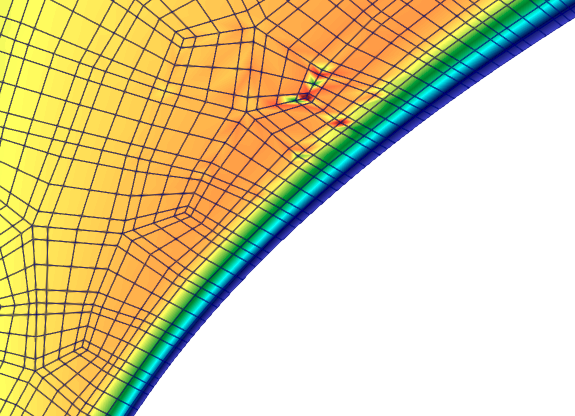}
        \caption{Pure Galerkin}
    \end{subfigure}
    \begin{subfigure}[b]{0.4\textwidth}
        \centering
        \includegraphics[width=\textwidth]{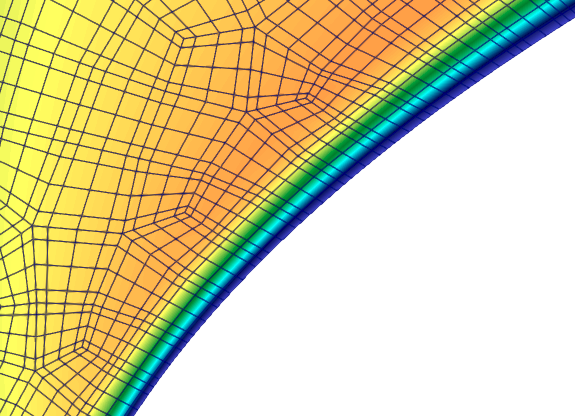}
        \caption{LPS}
    \end{subfigure}
    \begin{subfigure}[b]{0.06\textwidth}
        \centering
        \includegraphics[width=\textwidth]{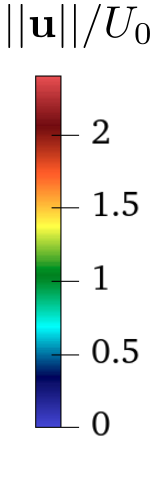}
    \end{subfigure}
    \caption{Velocity magnitude at the leading edge of the airfoil. Numerical instabilities are observed in the Pure Galerkin case.}
    \label{fig:naca_zoom}
\end{figure}

\begin{figure}[ht]
    \captionsetup[subfigure]{justification=centering}
    \centering
    \begin{subfigure}[b]{0.9\textwidth}
        \centering
        \includegraphics[width=\textwidth]{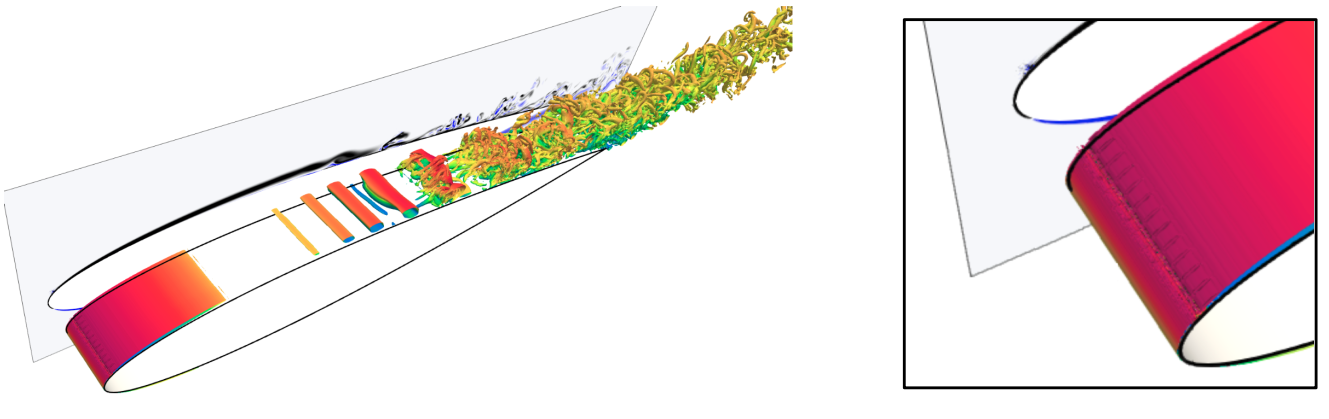}
        \caption{Pure Galerkin}
    \end{subfigure}
    \begin{subfigure}[b]{0.9\textwidth}
        \centering
        \includegraphics[width=\textwidth]{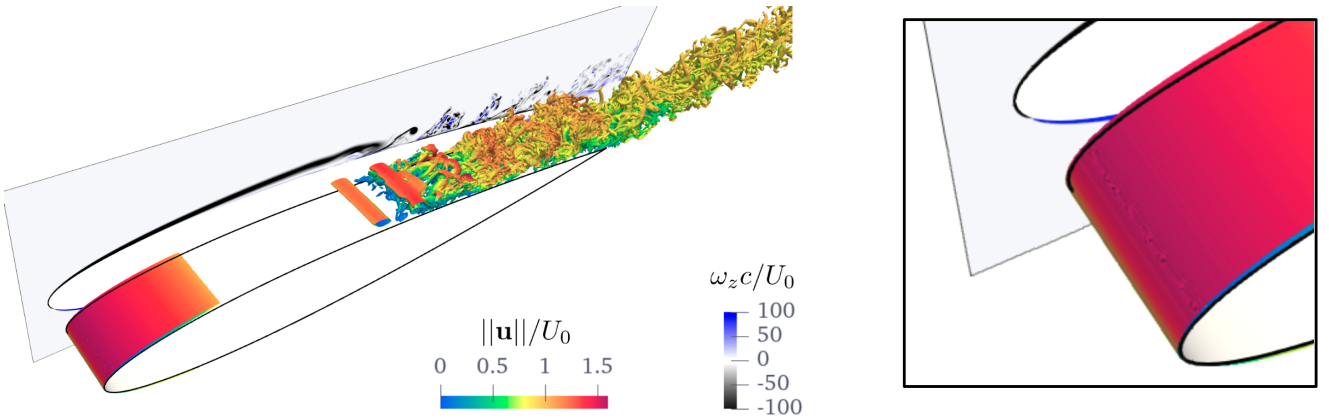}
        \caption{LPS}
    \end{subfigure}
    \caption{Q-vorticiy isocontour colored by the normalized velocity magnitude. The plane shows the normalized vorticity at the mid plane in spanwise direction.}
    \label{fig:naca_vorticity}
\end{figure}

The instantaneous solution is presented in Figure \ref{fig:naca_vorticity}, which shows the vorticity distribution around the airfoil.
The overall flow remains stable due to the inherent physical viscosity of the problem. However, as shown in Fig. \ref{fig:naca_zoom}, the Pure Galerkin method exhibits localized instabilities near the leading edge. In this region, the flow experiences strong acceleration and reaches the highest local Peclet number in the entire domain. At such high speed, the natural viscosity of the fluid is insufficient to fully stabilize the solution. Moreover, these instabilities are intensified by the irregularity of the mesh.
This is not the case of applying stabilization, which effectively eliminates the instabilities.

\begin{table}[ht]
\centering
\begin{tabular}{|c|c|c|c|c|}
\hline
Case & $C_l$ & $C_d$ & $x_\mathrm{sep}$ & $x_\mathrm{rea}$ \\
\hline
Lehmkuhl et al. \cite{Lehmkuhl2013} & 0.569 & 0.0291 & 0.169 & 0.566 \\
Zhang and Samtaney \cite{Zhang2016_airfoil} & 0.568 & 0.0285 & 0.141 & 0.580 \\
Present, Pure Galerkin  & 0.558 & 0.0241 & 0.176 & 0.517 \\
Present, LPS  & 0.558 & 0.0271 & 0.158 & 0.586 \\
\hline
\end{tabular}
\caption{Comparison of lift/drag coefficients and separation/reattachment distance in the streamline direction.}
\label{tab:naca}
\end{table}

According to Ref. \cite{Lehmkuhl2013}, the flow remains laminar up to the separation point, where it detaches from the surface and subsequently reattaches, forming a laminar separation bubble. Beyond this reattachment point, the flow transitions to a turbulent regime. The study of the separation mechanism and the correct prediction of the transition from laminar to turbulent flow in the detached shear-layer are both key aspects for improving engineering designs.

In the case of the Pure Galerkin method, instabilities at the leading edge appear to propagate downstream, disrupting the flow and altering the separation region. This results in multiple separation rollers and a delayed onset of separation, which leads to a higher drag coefficient compared to the results in the literature, as seen in Table \ref{tab:naca}.
Conversely, the stabilized case presents a proper capturing of the only laminar separation bubble and the transition to the turbulent regime starts earlier.

\begin{figure}[ht]
    \captionsetup[subfigure]{justification=centering}
    \centering
    \begin{subfigure}[b]{0.8\textwidth}
        \centering
        \includegraphics[width=\textwidth]{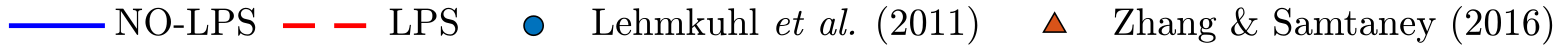}
    \end{subfigure}
    \begin{subfigure}[b]{0.49\textwidth}
        \centering
        \includegraphics[width=\textwidth]{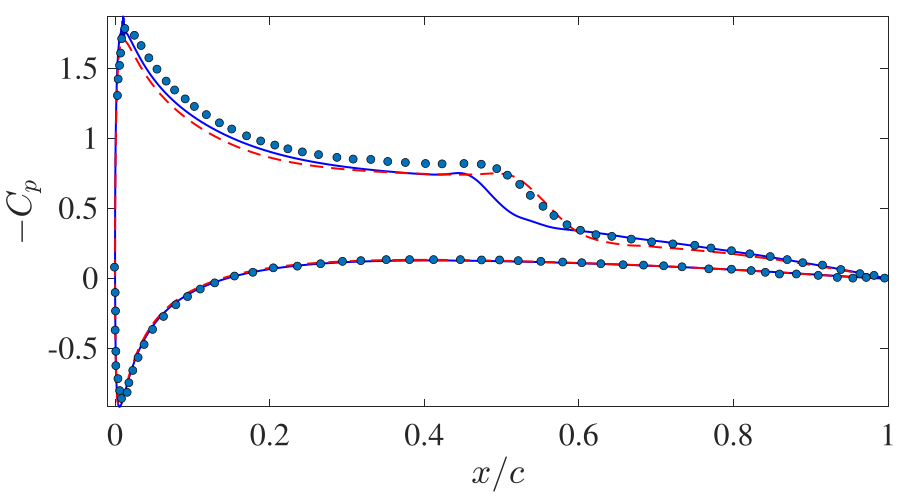}
        \caption{Pressure coefficient}
    \end{subfigure}
    \begin{subfigure}[b]{0.49\textwidth}
        \centering
        \includegraphics[width=\textwidth]{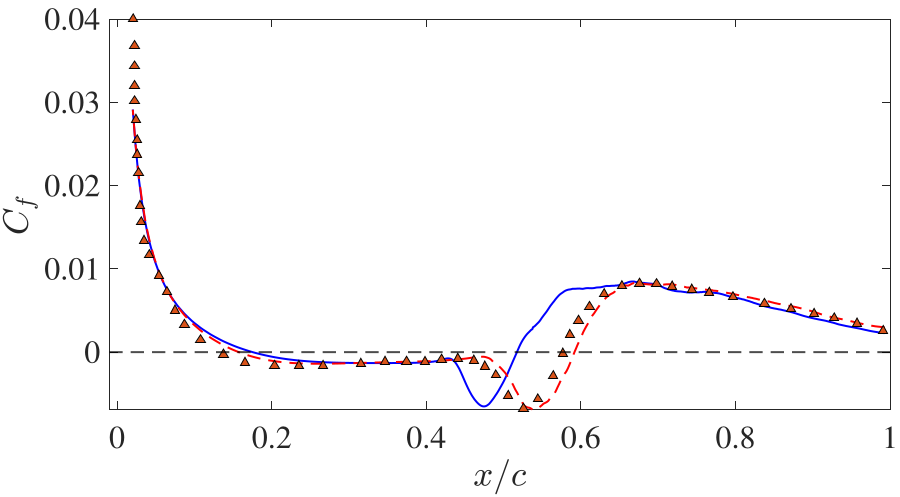}
        \caption{Friction coefficient}
    \end{subfigure}
    \caption{Distributions of the mean pressure coefficient and skin friction coefficient around the wall, plotted on the streamline direction normalized by the chord length.}
    \label{fig:naca_coefficients}
\end{figure}

Figure \ref{fig:naca_coefficients} shows the distribution of pressure and skin friction coefficients along the wall in the streamline direction. The case using LPS stabilization produces values that are in closer agreement with the reference solution, consistent with the results reported in Table \ref{tab:naca}. This highlights the effectiveness of the proposed stabilization method. Accurate prediction of the skin friction coefficient is particularly important, as it directly reflects boundary layer behavior and is critical for identifying flow separation, reattachment, and transition. Capturing these features reliably is essential for the validation and performance assessment of numerical methods in aerodynamic applications.

\section{Conclusions}

In this work, we have proposed a projection-based stabilization method based on the velocity field for solving incompressible flows. We demonstrate that the method exhibits convergence with both mesh refinement and increasing element order. It is also observed that the additional stabilization term eliminates aliasing of unresolved modes, enabling the use of the conservative form of the convection term. Nonuniform meshes, which are typically prone to sharp, non-physical oscillations between elements, benefit significantly from the proposed approach, which effectively suppresses these instabilities. While the stabilization introduces some numerical dissipation, this effect diminishes notably with higher-order elements, making the method especially well-suited for practical high-order spectral element simulations.

\section*{Acknowledgments}

The research leading to these results has received funding from the H2AERO project CPP2022-009921 from the Ministerio de Ciencia e Innovación.

V. Kumar acknowledges his AI4S fellowship within the Generaci\'on D initiative by Red.es, Ministerio para la Transformaci\'on Digital y de la Funci\'on P\'ublica, for talent attraction (C005/24-ED CV1), funded by NextGenerationEU through PRTR.

D. Mira acknowledges the grant Ayudas para contratos Ramón y Cajal (RYC 2021): RYC2021-034654-I from Ministerio de Ciencia e Innovación.

O. Lehmkuhl's work is financed by a Ramón y Cajal postdoctoral contract by the Ministerio de Economía y Competitividad, Secretaría de Estado de Investigación, Desarrollo e Innovación, Spain (RYC2018-025949-I). The authors acknowledge the support of the Departament de Recerca i Universitats de la Generalitat de Catalunya through the research group Large-scale Computational Fluid Dynamics (ref. 2021 SGR 00902).




\bibliographystyle{elsarticle-num} 
\bibliography{bibliography}



\end{document}